\documentclass[twocolumn]{autart}

\usepackage[english]{babel}
\usepackage{amsmath, amssymb}
\usepackage{savesym}
\savesymbol{AND}
\usepackage{algorithm}
\usepackage{algorithmic}
\usepackage{graphicx}
\usepackage[dvipsnames]{xcolor}
\usepackage{setspace}
\usepackage[hidelinks]{hyperref}
\usepackage{array,booktabs}
\usepackage{pifont}
\newcommand{\xmark}{\ding{55}}
\usepackage{flushend}

\allowdisplaybreaks

\theoremstyle{plain}
\newtheorem{ass}{Assumption}
\newtheorem{lemma}{Lemma}
\makeatletter
\newenvironment{proof}[0]{\par
  \normalfont \topsep6\p@\@plus6\p@\relax
  \trivlist
  \item[\hskip\labelsep
        \itshape
    Proof\@addpunct{.}]\ignorespaces
}{%
  \hfill\ \qed\endtrivlist\@endpefalse
}
\makeatother

\makeatletter
\newenvironment{proofof}[1]{\par
  \normalfont \topsep6\p@\@plus6\p@\relax
  \trivlist
  \item[\hskip\labelsep
        \bfseries
    Proof of #1\@addpunct{.}]\ignorespaces
}{%
  \hfill\ \qed\endtrivlist\@endpefalse
}
\makeatother

\newcommand{\bD}{\mathbb{D}}
\newcommand{\bN}{\mathbb{N}}
\newcommand{\bR}{\mathbb{R}}

\newcommand{\cD}{\mathcal{D}}
\newcommand{\cP}{\mathcal{P}}

\newcommand{\T}{\!^\top\!}
\newcommand{\x}{\times}

\newcommand{\norm}[1]{\| #1 \|}
\DeclareMathOperator{\dist}{dist}
\DeclareMathOperator{\relint}{relint}

\newcommand{\slater}{\tilde{x}}
\newcommand{\xh}{\hat{x}}
\newcommand{\xs}{x^\star}
\newcommand{\ls}{\lambda^\star}

\newcommand{\Xs}{X^\star}
\newcommand{\Ls}{\Lambda^\star}

\newcommand{\xhi}{\xh_i}
\newcommand{\li}{\lambda_i}
\newcommand{\lj}{\lambda_j}
\newcommand{\ai}{\ell_i}

\newcommand{\vk}{v(k)}
\newcommand{\xhk}{\xh(k)}

\newcommand{\xhik}{\xhi(k)}
\newcommand{\eik}{e_i(k)}
\newcommand{\lik}{\li(k)}
\newcommand{\ljk}{\lj(k)}
\newcommand{\aik}{\ai(k)}

\newcommand{\vkp}{v(k+1)}
\newcommand{\xkp}{x(k+1)}
\newcommand{\xikp}{x_i(k+1)}
\newcommand{\xhkp}{\xh(k+1)}
\newcommand{\xhikp}{\xhi(k+1)}
\newcommand{\eikp}{e_i(k+1)}

\newcommand{\likp}{\li(k+1)}

\newcommand{\vrp}{v(r+1)}
\newcommand{\xrp}{x(r+1)}
\newcommand{\lir}{\li(r)}
\newcommand{\ljr}{\lj(r)}
\newcommand{\air}{\ai(r)}
\newcommand{\xirp}{x_i(r+1)}
\newcommand{\lirp}{\li(r+1)}

\begin{document}

\begin{frontmatter}
\runtitle{Dual decomposition for multi-agent distributed optimization with coupling constraints}  

\title{Dual decomposition for multi-agent\\ distributed optimization with coupling constraints}

\thanks{Research was supported by the European Commission under the project UnCoVerCPS, grant number 643921.}

\author[PoliMi]{Alessandro~Falsone}\ead{\\ alessandro.falsone@polimi.it},
\author[Oxford]{Kostas~Margellos}\ead{kostas.margellos@eng.ox.ac.uk},
\author[PoliMi]{Simone~Garatti}\ead{simone.garatti@polimi.it},
\author[PoliMi]{Maria~Prandini}\ead{maria.prandini@polimi.it}

\address[PoliMi]{Dipartimento di Elettronica, Informazione e Bioingegneria, Politecnico di Milano, Via Ponzio 34/5, 20133 Milano, Italy}
\address[Oxford]{Department of Engineering Science, University of Oxford, Parks Road, Oxford, OX1 3PJ, United Kingdom}


\begin{abstract}
	We study distributed optimization in a cooperative multi-agent setting, where agents have to agree on the usage of shared resources and can communicate via a time-varying network to this purpose. Each agent has its own decision variables that should be set so as to minimize its individual objective function subject to local constraints. Resource sharing is modeled via coupling constraints that involve the non-positivity of the sum of agents' individual functions, each one depending on the decision variables of one single agent.
	We propose a novel distributed algorithm to minimize the sum of the agents' objective functions subject to both local and coupling constraints, where dual decomposition and proximal minimization are combined in an iterative scheme. Notably, privacy of information is guaranteed since only the dual optimization variables associated with the coupling constraints are exchanged by the agents.
	Under convexity assumptions, jointly with suitable connectivity properties of the communication network, we are able to prove that agents reach consensus to some optimal solution of the centralized dual problem counterpart, while primal variables converge to the set of optimizers of the centralized primal problem. The efficacy of the proposed approach is demonstrated on a plug-in electric vehicles charging problem.
\end{abstract}

\end{frontmatter}

\section{Introduction}
This paper addresses optimization in multi-agent networks where each agent aims at optimizing a local performance criterion possibly subject to local constraints, but yet needs to agree with the other agents in the network on the value of some decision variables that refer to the usage of some shared resources.

Cooperative multi-agent decision making problems have been studied recently by many researchers, mainly within the control and operational research communities, and are found in various application domains such as  power systems \cite{bolognani_etal_2015,zhang_giannakis_2015}, wireless and social networks \cite{mateos_giannakis_2012,baingana_etal_2014}, robotics \cite{martinez_etal_2007}, to name a few.

A possible approach to cooperative multi-agent optimization consists in formulating and solving a mathematical program involving the decision variables, objective functions, and constraints of the entire network. Though this centralized perspective appears sensible, it may end up being impractical for large scale systems for which the computational effort involved in the program solution can be prohibitive. Also, privacy of information is not preserved since agents are required either to share among them or to provide to a central entity their performance criteria and constraints.

Distributed optimization represents a valid alternative to centralized optimization and, in particular, it overcomes the above limitations by allowing agents to keep their information private, while distributing the computational effort. Typically, an iterative procedure is conceived, where at each iteration agents perform some local computation based on their own information and on the outcome of the local computations of their neighboring agents at the previous iteration, till convergence to some solution, possibly an optimal one for the centralized optimization problem counterpart.

Effective  distributed optimization algorithms have been proposed in the literature for a general class of convex problems over time-varying, multi-agent networks. In particular, consensus-based optimization algorithms are formulated in \cite{Nedic_Ozdaglar_2009,nedic_etal_2010,Ram_etal_2012,Lee_Nedic_2013} and in our recent paper \cite{our_paper} for problems where agents have their own objective functions and constraints but decision variables are common.

In this paper, we address a specific class of convex optimization problems over time-varying, multi-agent networks, which we refer to as \emph{inequality-coupled problems} for short-hand notation. In this class of problems, each agent has its own decision vector, objective function, and constraint set, and is coupled to the others via a constraint expressed as the non-positivity of the sum of convex functions, each function corresponding to one agent.

We propose a novel distributed iterative scheme based on a combination of dual decomposition and proximal minimization to deal with inequality-coupled problems.
Under convexity assumptions and suitable connectivity properties of the communication network, agents reach
consensus with respect to the dual variables, without disclosing information about their optimal decision, local objective and constraint functions, nor about the function encoding their contribution to the coupling constraint.
The proposed algorithm converges to some optimal dual solution of the centralized problem counterpart, while for the primal variables, we show convergence to the set of optimal primal solutions.

The contributions of our paper versus the existing literature are summarized in the following.

\begin{table*}[t]
	\centering
	\begin{tabular}{m{22em}ccccc}
		\toprule
			& \cite{our_paper}
			& \cite{zhu_martinez_2012}
			& \cite{nedic_scaglione_2014}
			& \cite{simonetto2016primal}
			& Algorithm~\ref{alg:Alg1}		\\ \midrule
		Convergence to optimal solution
			& \checkmark
			& \checkmark
			& \checkmark
			& \xmark
			& \checkmark	\\
		Exchange of info related to
			& P
			& P+D
			& P+D
			& D
			& D				\\
		Required technical assumptions
			& CC
			& CC+SP
			& CC+SP+DF
			& CC+SP
			& CC			\\
		Type of coupling constraints
			& $\leq$
			& $\leq$
			& $\leq$
			& $\leq$
			& $\leq+=$		\\
		Network topology
			& V
			& V
			& V
			& F
			& V				\\
		Number of decision variables in the local problem
			& $\sum_i n_i$
			& $\sum_i n_i$
			& $n_i$
			& $n_i$
			& $n_i$			\\
		Number of variables that need to be stored locally 
			& $\sum_i n_i$
			& $\sum_i n_i + p$
			& $3(n_i+p)$
			& $2n_i + p$
			& $2n_i + p$	\\
		Number of variables that are exchanged within any communicating agents pair 
			& $\sum_i n_i$
			& $\sum_i n_i + p$
			& $2p$
			& $p$
			& $p$			\\ \bottomrule
	\end{tabular}
	\caption{Comparison against other approaches. Legend: + means ``and'', P stands for primal, D stands for dual, CC stands for convexity and compactness, SP stands for knowledge of a Slater point by agents, DF stands for differentiability of the objective function, $\leq$ stands for inequality constraints, $=$ stands for equality constraints, V stands for time-varying, F stands for fixed, $n_i$ is the number of decision variables of agent $i$, and $p$ is the number of coupling constraints/Lagrange multipliers.}
	\label{tab:alg_comparison}
\end{table*}

Our scheme can be seen as an extension of dual decomposition based algorithms to a distributed setting, accounting for time-varying network connectivity. As a matter of fact, if the communication networks were time-invariant and connected, then, dual decomposition techniques (see \cite{Bo_Johansson_2010}, and references therein) as well as approaches based on the alternating direction method of multipliers \cite{Boyd_etal_2010,shi_ling_2014} could be applied to the set-up of this paper, since, after dualizing the coupling constraint, the problem assumes a separable structure. However, in \cite{Bo_Johansson_2010} and \cite{Boyd_etal_2010} a central update step involving communication among all agents that are coupled via the constraints is required for the dual variables, and this prevents their usage in the distributed with time-varying connectivity set-up of this paper. In \cite{shi_ling_2014} no central update step is needed but the constraints appearing in the dual problem cannot be handled. An interesting distributed dual decomposition based algorithm which overcomes the need for a central node and which is more in line with our scheme has been proposed in \cite{simonetto2016primal}. The main differences between \cite{simonetto2016primal} and our algorithm are as follows:
\begin{itemize}
	\item[a.] the algorithm of \cite{simonetto2016primal} requires that the communication network is time invariant, while our algorithm admits time-variability;
	\item[b.] in \cite{simonetto2016primal} a constant step-size is employed, while our algorithm uses a vanishing step-size. The constant step-size has the advantage of enhancing a faster convergence rate, but, at the same time, convergence to a neighborhood of the optimal is guaranteed only. Our algorithm instead is guaranteed to converge to the optimal solution of the original problem;
	\item[c.] the algorithm of \cite{simonetto2016primal} requires that a Slater point exists and is known to all agents, while existence only is required in our algorithm. This relaxation of the conditions for the applicability of the approach can be crucial in those cases where a Slater point is not a-priori available since the reconstruction of a Slater point in a distributed set-up seems to be as challenging as the original problem and requires extra synchronization among agents.
\end{itemize}

From another perspective, which is better explained later on in the paper, our approach can be also interpreted as a subgradient based algorithm for the resolution of the dual problem, equipped with an auxiliary sequences that allows one to recover the solution of the primal problem we are interested in.
In this respect related contributions are \cite{Nedic_Bertsekas_2001,bertsekas_2003,Bertsekas_2011}, where some incremental gradient/subgradient algorithms that can be adopted as an alternative to dual decomposition are proposed. These algorithms, however, require that agents perform updates sequentially, in a cyclic or randomized order, and do not really fit the distributed set-up of this paper. The recent contributions \cite{zhu_martinez_2012} and \cite{nedic_scaglione_2014} instead present primal-dual subgradient based consensus algorithms that fit our set-up and are comparable to our approach. The main differences are:
\begin{itemize}
	\item[d.] in \cite{zhu_martinez_2012} a global knowledge by all agents of the coupling constraint in the primal is required and in both \cite{zhu_martinez_2012} and \cite{nedic_scaglione_2014} information related to the primal problem is exchanged among agents while the algorithm is running. In the separable set-up of this paper, the agents local information on the primal problem (namely, the value of the local optimization variables, the local objective function, the local constraints, and the contribution of the agent to the coupling constraint) can be regarded as sensitive data and their exchange as in \cite{zhu_martinez_2012} and \cite{nedic_scaglione_2014} may raise privacy issues. In our algorithm, only the local estimates of the dual variables are exchanged, and this secures maximum privacy among agents;
	\item[e.] the algorithms of \cite{zhu_martinez_2012} and \cite{nedic_scaglione_2014} require that a Slater point exists and is known to all agents, while existence only is required in our algorithm. As commented before, requiring the knowledge of a Slater point by the agents can hamper the usability of the algorithm. Moreover, the convergence to the optimal solution in \cite{nedic_scaglione_2014} is guaranteed only when each agents objective function is differentiable;
	\item[f.] to apply the algorithm of \cite{zhu_martinez_2012} to our set-up, each agent has to generate local copies of the optimization variables of all the other agents, which then are optimized and exchanged. This often results in an unnecessary increase of the computational and communication efforts, which indeed scale as the number of agents in the network. In our approach instead agents need to optimize the local variables only and exchange the estimate of the dual variables, which are as many as the number of coupling constraints. The required local computational effort is thus much smaller. As for the communication effort, our approach is particularly appealing when the number of coupling constraints is low compared to the overall dimensionality of primal decision variables.
\end{itemize}

Finally, note that the approaches to distributed optimization in \cite{Nedic_Ozdaglar_2009,nedic_etal_2010,Ram_etal_2012,Lee_Nedic_2013,our_paper}, which do not resort to any dual problem, can be applied to inequality-coupled problems by introducing  a common decision vector collecting all local decision variables. This, however, immediately leads to the drawback of an increased computational and communication effort as discussed in point f above. Moreover, these approaches requires an exchange of information related to the primal, which leads to the privacy issues outlined in point d above.

Table~\ref{tab:alg_comparison} summarizes the comparison between the proposed methodology and the most significant approaches that apply to the same set-up. In the table, algorithms are assessed each against the others based on several indices related to points a-f above. Perhaps, it is worth mentioning that, since in \cite{our_paper} and \cite{zhu_martinez_2012} local copies of the optimization variables of all agents are required, a further issue arises for these two algorithms. As matter of fact, since agent $i$ has no constraints for the variables of the other agents, the assumption, which is common to all algorithms, of compactness of the overall optimization domain is no longer verified. This issue can be prevented e.g. by forcing $x_j$ to belong to an outer box to the constraints set $X_j$ for each $j\neq i$, but in doing so some information about the local $X_i$ is exchanged, leading to further privacy issues.

A preliminary version of this work is given in \cite{falsone2016cdc}. The present paper significantly extends that contribution from a theoretical viewpoint, in that it contains the proofs of all the results stated in the conference version and a preliminary study on the convergence rate. Furthermore, a thorough comparison with the literature has been added and an assessment of the performance of the proposed approach has been carried out through a concrete problem on plug-in electric vehicles charging, where also guidelines on how to speed up numerical convergence are provided.

The rest of the paper is organized as follows. A formal statement of the problem, along with the proposed algorithm is given in Section~\ref{sec:problem}. Convergence and optimality of our algorithm are studied in Section~\ref{sec:analysis}. A numerical example along with a suggestion on how to speed up numerical convergence is presented in Section~\ref{sec:example}, while some final remarks are drawn in Section~\ref{sec:conclusion}.

\section{Distributed constrained optimization} \label{sec:problem}
\subsection{Problem statement and proposed solution} \label{sec:setup}
Consider the following optimization program
\begin{equation}\label{eq:cP}
	\begin{aligned}
		\cP: \;\; \min_{{\{x_i\in X_i\}}_{i=1}^m} \quad &\sum_{i=1}^m f_i(x_i) \\
			\text{subject to:} \quad &\sum_{i=1}^m g_i(x_i) \leq 0,
	\end{aligned}
\end{equation}
involving  $m$ agents that communicate over a time-varying network.
Each agent $i$, $i=1,2,\dots,m$, has its own vector $x_i\in\bR^{n_i}$ of $n_i$ decision variables, its  local constraint set $X_i\subseteq \bR^{n_i}$ and  objective function $f_i(\cdot):\bR^{n_i}\rightarrow\bR$, and it is contributing to the coupling constraint $\sum_{i=1}^m g_i(x_i) \leq 0$ via function $g_i(\cdot):\bR^{n_i}\rightarrow \bR^p$. Note that equality linear coupling constraints can be also dealt with by means of $\cP$, by means of double-sided inequalities.

Problem $\cP$ could be solved, in principle, in a centralized fashion. However, if the number $m$ of agents  is large, this may turn out to be computationally prohibitive. In addition, each agent would be required to share its own information (coded via $f_i(\cdot)$, $X_i$, and $g_i(\cdot)$) either with the other agents or with a central unit collecting all information, which may be undesirable in some cases, due to privacy issues.

We next formulate a distributed strategy that overcomes both the privacy and computational issues outlined above by resorting to the dual of \eqref{eq:cP}.

Let us consider the Lagrangian function $L(x,\lambda): \mathbb{R}^n \times \mathbb{R}^p_+ \to \mathbb{R}$ given by
\begin{equation*}
	L(x,\lambda) = \sum_{i=1}^m L_i(x_i,\lambda) = \sum_{i=1}^m \left\{ f_i(x_i) + \lambda\T g_i(x_i) \right\},
\end{equation*}
where $x=[x_1\T\,\cdots\,x_m\T\,]\,\T\in X = X_1\times\cdots\times X_m \subseteq\bR^n$, with $n=\sum_{i=1}^m n_i$, whereas $\lambda\in\bR^p_+$ is the vector of Lagrange multipliers ($\bR^p_+$ denotes the $p$-th dimensional non-negative orthant; in the sequel we shall sometimes write $\lambda\geq0$ in place of $\lambda\in\bR^p_+$).

Correspondingly, we can define the dual function as
\begin{equation}
	\varphi(\lambda) = \min_{x\in X} L(x,\lambda), \label{eq:dual_function}
\end{equation}
which, due to the separable structure of objective and constraint functions in problem $\cP$ (see \eqref{eq:cP}), can be expressed as
\begin{equation}\label{eq:dual_i}
	\varphi(\lambda) = \sum_{i=1}^m \varphi_i(\lambda) = \sum_{i=1}^m \min_{x_i\in X_i} L_i(x_i,\lambda),
\end{equation}
where each $\varphi_i(\cdot)$ is a concave function representing the dual function of agent $i$.

Given these definitions, the dual of problem $\cP$ in \eqref{eq:cP} can be expressed as:
\begin{equation*}
	\cD: \;\; \max_{\lambda \geq 0} \min_{x\in X} L(x,\lambda),
\end{equation*}
or, equivalently, as
\begin{equation}
	\cD: \;\; \max_{\lambda \geq 0} \sum_{i=1}^m \varphi_i(\lambda). \label{eq:cD}
\end{equation}
The coupling between agents is given in \eqref{eq:cD} by the fact that $\lambda$ is a common decision vector and the agents should agree on its value.

\begin{algorithm}[t]
	\caption{Distributed algorithm}
	\begin{spacing}{1.1}
	\begin{algorithmic}[1]
		\STATE \textbf{Initialization} \\
		\STATE ~~~~$k=0$. \\
        \STATE ~~~~Consider $\xhi(0)\in X_i$, for all $i=1,\dots,m$. \\
		\STATE ~~~~Consider $\li(0)\in\bR^p_+$, for all $i=1,\dots,m$. \\
		\STATE \textbf{For $i=1,\ldots,m$ repeat until convergence} \\
		\STATE ~~~~$\aik = \sum_{j=1}^m a_j^i(k) \ljk$. \\
		\STATE ~~~~$\xikp \in \arg\min_{x_i \in X_i} f_i(x_i) + \aik\T g_i(x_i)$. \\
		\STATE ~~~~$\likp = \arg\max_{\lambda_i \geq 0} \big\{ g_i(\xikp)\T \lambda_i$ \\
			\begin{flushright} $- \frac{1}{2c(k)}\norm{\lambda_i - \aik}^2 \big\}$ \end{flushright}
		\STATE ~~~~$\xhikp = \xhik + \frac{c(k)}{\sum_{r=0}^k c(r)} (\xikp-\xhik)$.
		\STATE ~~~~$k \gets k+1$.
	\end{algorithmic}
	\end{spacing}
	\label{alg:Alg1}
\end{algorithm}

%

Algorithm~\ref{alg:Alg1} is a distributed iterative scheme that aims at reconstructing the solution to both the dual problem \eqref{eq:cD} and the primal problem \eqref{eq:cP} by exchanging a minimal amount of information among agents. Its steps are explained hereafter.

Each agent $i$, $i=1,\dots,m$, initializes the estimate of its local decision vector with $\xhi(0)\in X_i$ (step~3 of Algorithm~\ref{alg:Alg1}), and the estimate of the common dual variables vector with a $\li(0)\in\bR^p_+$ that is feasible for problem $\cD$ (step~4 of Algorithm~\ref{alg:Alg1}).
A sensible choice is to set $\xhi(0) \in \arg \min_{x_i \in X_i} f_i(x_i)$, and $\li(0) = 0$, $i=1,\ldots,m$, which corresponds to the solution of problem \eqref{eq:cP} when coupling constraints are neglected.

At every iteration $k$, $k\ge 1$, each agent $i$ computes a weighted average $\aik$ of the dual variables vector based on the estimates $\ljk$, $j=1,\dots,m$, of the other agents and its own estimate (step~6). The weight $a_j^i(k)$ that agent $i$ attributes to the estimate of agent $j$ at iteration $k$ is set equal to zero if agent $i$ does not communicate with agent $j$ at iteration $k$.

Then, Algorithm~\ref{alg:Alg1} alternates between a primal and a dual update step (step~7 and step~8, respectively). In particular, step~7 performs an update of the local primal vector $\xikp$ by minimizing $L_i$ evaluated at $\lambda = \aik$ as in dual decomposition, whereas, differently from dual decomposition, which would consists of the maximization of $L_i$ evaluated at $x_i = x_i(k+1)$, the update of the dual vector in step~8  involves also the proximal term $- \frac{1}{2c(k)}\norm{\lambda_i - \aik}^2$ to foster consensus among the agents.

Steps~7 and~8 can be thought of as an approximation of the following proximal maximization step
\begin{align} \label{eq:full_proxi}
	\likp = \arg\max_{\lambda_i \geq 0} \min_{x_i\in X_i} \Big\{ L_i&(x_i,\lambda_i) \\
		&- \frac{1}{2c(k)}\norm{\lambda_i - \aik}^2 \Big\}, \nonumber
\end{align}
which would implement the distributed algorithm of \cite{our_paper} for the dual problem \eqref{eq:cD}. Steps~7 and~8 are however preferred to \eqref{eq:full_proxi} since the resolution of the $\max-\min$ program in \eqref{eq:cD} is very hard in general. Moreover, it is perhaps worth mentioning at the outset that step~8 in Algorithm~\ref{alg:Alg1} is equivalent to a projected subgradient step. Indeed the constrained maximization of a quadratic function in step~8 can be explicitly solved, leading to
\begin{equation} \label{eq:subgradient_step}
	\likp = [\aik + c(k) g_i(\xikp)]_+,
\end{equation}
where $[\,\cdot\,]_+$ denotes the projection of its argument onto $\bR^p_+$.
Then, it can be shown that $g_i(\xikp)$ is a subgradient of the dual function $\varphi_i(\cdot)$ evaluated at $\aik$ (see the proof of Theorem \ref{thm:dual_optimality} for more details), while $c(k)$ can be thought of as the subgradient step-size. Hence, steps~7 and~8 can be also seen as an application of the distributed subgradient algorithm of \cite{nedic_etal_2010}, which was originally developed for primal problems though, to the dual problem \eqref{eq:cD}.


Unfortunately, the local primal vector $x_i(k)$ does not converge to the optimal solution $x^\ast_i$ to \eqref{eq:cP} in general. Therefore, the auxiliary primal iterates $\xhikp$, defined as the weighted running average of $\{\xirp\}_{r=0}^k$
\begin{equation} \label{eq:xhat_def}
	\xhikp = \frac{\sum_{r=0}^k c(r) \xirp}{\sum_{r=0}^k c(r)},
\end{equation}
is computed in step 9 of Algorithm~\ref{alg:Alg1} in a recursive fashion. Such an auxiliary variable shows better convergence properties as compared to $x_i(k)$, and is often constructed in the so-called primal recovery procedure of dual decomposition methods, \cite{nedic_ozdaglar_2009_approx,nedic_scaglione_2014,zhu_martinez_2012}.

Note that in Algorithm \ref{alg:Alg1} no local information related to the primal is exchanged between the agents (as a matter of fact only the estimates of the dual vector are communicated) so that our algorithm is well suited to account for privacy requirements.


\subsection{Structural and communication assumptions} \label{sec:assumptions}
The proposed distributed algorithm shows properties of convergence and optimality, which hold under the following assumptions on the structure of the problem and on the communication features of the time-varying multi-agent network.
\begin{ass}{[Convexity]} \label{ass:convex}
	For each $i=1,\dots,m$, the function $f_i(\cdot):\bR^{n_i}\to\bR$ and each component of $g_i(\cdot):\bR^{n_i} \to \bR^p$ are convex; for each $i=1,\dots,m$ the set $X_i\subseteq\bR^{n_i}$ is convex.
\end{ass}

\begin{ass}{[Compactness]} \label{ass:compact}
	For each $i=1,\ldots,m$, the set $X_i\subseteq\bR^{n_i}$ is compact.
\end{ass}
Note that, under Assumptions~\ref{ass:convex} and \ref{ass:compact}, $\norm{g_i(x_i)}$ is finite for any $x_i\in X_i$: $\norm{g_i(x_i)} \leq G$, $\forall x_i\in X_i$, where $G = \max_{i=1,\dots,m} \max_{x_i\in X_i} \norm{g_i(x_i)}$.

\begin{ass}{[Slater's condition]} \label{ass:slater}
There exists $\slater=[\slater_1\,\cdots\,\slater_m]\T\in \relint(X)$, where $\relint(X)$ is the relative interior of the set $X$, such that $\sum_{i=1}^m g_i(\slater_i) \leq 0$ for those components of $\sum_{i=1}^m g_i(x_i)$ that are linear in $x$, if any, while $\sum_{i=1}^m g_i(\slater_i) < 0$ for all other components.
\end{ass}
As a consequence of Assumptions~\ref{ass:convex}-\ref{ass:slater}, we have that strong duality holds and an optimal primal-dual pair $(\xs,\ls)$ exists, where $\xs=[\xs_1\,\cdots\,\xs_m]\T$. Moreover, the Saddle-Point Theorem holds, \cite{boyd_2004}, i.e., given an optimal pair $(\xs,\ls)$, we have that
\begin{equation} \label{eq:saddle_point}
	L(\xs,\lambda) \leq L(\xs,\ls) \leq L(x,\ls), \ \lambda\in\bR^p_+, \ x\in X.
\end{equation}

The reader should note that, differently from other approaches, we require a Slater point to exists, but we do not need the agents to actually compute it, which, as discussed in the introduction, might be impractical in a distributed set-up.

In the following we will denote by $\Xs$ the set of all primal minimizers, and by $\Ls$ the set of all dual maximizers.

As for the time-varying coefficient $c(k)$, we impose the following assumptions that are similar to those in \cite{nedic_etal_2010,zhu_martinez_2012,our_paper}.
\begin{ass}{[Coefficient $c(k)$]} \label{ass:ck_coefficient}
$\{c(k)\}_{k \geq 0}$ is a non-increasing sequence of positive reals such that $c(k) \leq c(r)$ for all $k \geq r$, with $r \geq 0$. Moreover,
	\begin{enumerate}
		\item $\sum_{k=0}^{\infty} c(k) = \infty$,
		\item $\sum_{k=0}^{\infty} c(k)^2 < \infty$.
	\end{enumerate}
\end{ass}
One possible choice for $\{c(k)\}_{k \geq 0}$ satisfying Assumption~\ref{ass:ck_coefficient} is $c(k) = \beta/(k+1)$ for some $\beta>0$.

As in \cite{tsitsiklis_1984,tsitsiklis_etal_1986,olshevsky_tsitsiklis_2011}, the communication network is required to satisfy the following connectivity conditions.
\begin{ass}{[Weight coefficients]} \label{ass:weights}
	There exists $\eta \in (0,1)$ such that for all $i,j \in \{1,\dots,m\}$ and all $k \geq 0$, $a_j^i(k) \in [0,1)$, $a_i^i(k) \geq \eta$, and $a_j^i(k) > 0$ implies that $a_j^i(k) \geq \eta$.
	Moreover, for all $k \geq 0$,
	\begin{enumerate}
		\item $\sum_{j=1}^{m} a_j^i(k) = 1$ for all $i=1,\ldots,m$,
		\item $\sum_{i=1}^{m} a_j^i(k) = 1$ for all $j=1,\ldots,m$.
	\end{enumerate}
\end{ass}
Note that, if we fix $k \geq 0$, the information exchange between the $m$ agents can be coded via a directed graph $(V,E_k)$, where nodes in $V = \{1,\ldots,m\}$ represent the agents, and the set $E_k$ of directed edges is defined as
\begin{align} \label{eq:edges}
	E_k = \big \{ (j,i):~ a_j^i(k) > 0 \big \},
\end{align}
i.e., at time $k$ the link $(j,i)$ is present if agent $j$ communicates with agent $i$ and agent $i$ weights the information received from agent $j$ with $a_j^i(k)$. If the communication link is not active, then $a_j^i(k) = 0$; if $a_j^i(k) > 0$ then agent $j$ is said to be neighbor of agent $i$ at time $k$.

Let $E_{\infty} = \big \{ (j,i):~ (j,i) \in E_k \text{ for infinitely many } k \big \}$ denote the set of edges $(j,i)$ representing pairs of agents that communicate directly infinitely often. We then impose the following connectivity and communication assumption.
\begin{ass}{[Connectivity and communication]} \label{ass:network}
	Graph $(V,E_{\infty})$ is strongly connected, i.e., for any two nodes there exists a path of directed edges that connects them. Moreover, there exists $T \geq 1$ such that for every $(j,i) \in E_{\infty}$, agent $i$ receives information from a neighboring agent $j$ at least once every consecutive $T$ iterations.
\end{ass}

Details on the interpretation of Assumptions~\ref{ass:weights} and \ref{ass:network} can be found in \cite{Nedic_Ozdaglar_2009,our_paper,nedic_etal_2010}.

\subsection{Statement of the main results} \label{sec:solution}

If Assumptions \ref{ass:convex}-\ref{ass:network} are satisfied, then Algorithm \ref{alg:Alg1} converges and agents agree to a common vector of Lagrange multipliers. Specifically, their local estimates $\lik$ converge to some optimal vector of Lagrange multipliers, while the vector $\xhk = [\hat{x}_1(k)\T\,\cdots\,\hat{x}_m(k)\T]\T$ approaches $\Xs$, the set of minimizers of the primal problem.

These results are formally stated in the following theorems.
\begin{thm}{[Dual Optimality]} \label{thm:dual_optimality}
	Under Assumptions~\ref{ass:convex}-\ref{ass:network}, there exists a $\ls\in\Ls$ such that
	\begin{equation}
		\lim_{k \rightarrow \infty} \norm{\lik - \ls} = 0, \text{ for all } i=1,\dots,m. \label{eq:dual_optimality}
	\end{equation}
\end{thm}

\begin{thm}{[Primal Optimality]} \label{thm:primal_optimality}
	Under Assumptions~\ref{ass:convex}-\ref{ass:network}, we have that
	\begin{equation}
		\lim_{k \rightarrow \infty} \dist(\xhk,\Xs) = 0, 
	\end{equation}
	where $\dist(y,Z)$ denotes the distance between $y$ and the set $Z$, i.e., $\dist(y,Z) = \min_{z\in Z} \norm{y-z}$.
\end{thm}

\section{Convergence and Optimality Analysis} \label{sec:analysis}
This section is devoted to the convergence and optimality analysis of Algorithm~\ref{alg:Alg1}. We will first prove Theorem~\ref{thm:dual_optimality} employing the convergence result of the primal algorithm proposed in \cite{nedic_etal_2010} applied to \eqref{eq:cD}. We will then provide some preliminary results which are instrumental for the proof of Theorem~\ref{thm:primal_optimality}, and finally we will give the proof of Theorem~\ref{thm:primal_optimality}.

\subsection{Proof of Theorem~\ref{thm:dual_optimality}} \label{sec:proof_of_thm1}

The structure of problem \eqref{eq:cD} fits the framework considered in \cite{nedic_etal_2010}, and as already noted below equation \eqref{eq:subgradient_step}, the part of Algorithm \ref{alg:Alg1} that pertains to the update of the dual vector (namely, steps 6 and 7 and step 8 which is equivalent to \eqref{eq:subgradient_step}) is an implementation of the subgradient algorithm of \cite{nedic_etal_2010} for the dual problem \eqref{eq:cD}. In particular, referring to \eqref{eq:subgradient_step}, the fact that $g_i(\xikp)$, with $\xikp$ computed as in step~7, is a subgradient of $\varphi_i(\lambda) = \min_{x_i\in X_i} \{ f_i(x_i) + \lambda\T g_i(x_i) \}$ evaluated at $\lambda = \aik$ is a well-known consequence of the Danskin's theorem (see Proposition~B.25 in \cite{bertsekas_1999}). Moreover, since $f_i$ are convex over the whole domain $\bR^{n_i}$ by Assumption \ref{ass:convex} and since $\xikp\in X_i$, which is compact by Assumption~\ref{ass:compact}, it holds that the subgradients of each agent objective function evaluated at $\xikp$ are always bounded. This latter observation along with Assumptions~\ref{ass:convex}-\ref{ass:network} allows one to conclude that all requirements for Proposition~4 in \cite{nedic_etal_2010} to hold are verified, and then the result \eqref{eq:dual_optimality} of Theorem~\ref{thm:dual_optimality} follows by a direct application of Proposition~4 in \cite{nedic_etal_2010}. This concludes the proof.\hfill\ \qed

By Theorem~\ref{thm:dual_optimality}, for all $i=1,\dots,m$, the sequence $\{ \lik \}_{k \geq 0}$ is converging to some $\ls\in\Ls$. Therefore,  $\{ \lik \}_{k \geq 0}$ is also a bounded sequence, that is $\norm{\lik} \leq D$, with $D = \max_{i}\sup_{k\geq 0}\norm{\lik} < \infty$.

\subsection{Error relations} \label{sec:prep_results}
In this subsection we prove some preliminary relations that link the dual variables local estimate $\lik$ of agent $i$, $i=1,\ldots,m$, with their arithmetic average
\begin{equation}
	\vk = \frac{1}{m} \sum_{i=1}^m \lik, \text{ for all } k \geq 0. \label{eq:vk}
\end{equation}
These relations are then used in Section~\ref{sec:proof_of_thm2} for the proof of Theorem~\ref{thm:primal_optimality}.

To start with, the following lemma establishes a link between $\norm{\likp - \vkp}$ and $\norm{\eikp}$, $i=1,\ldots,m$, where
\begin{equation} \label{eq:def_error}
	\eikp = \likp-\aik
\end{equation}
is the consensus error for agent $i$.
\begin{lemma} \label{lemma:sum_error}
	Consider Assumptions~\ref{ass:ck_coefficient}-\ref{ass:network}. Fix any $\alpha_1\in\bR_+\setminus\{0\}$. We then have that for any $N \in \bN_+ \setminus \{0\}$,
	\begin{align}
		2 &\sum_{k=1}^N c(k) \sum_{i=1}^m \norm{\likp - \vkp} \nonumber  \\
		& < \alpha_1 \sum_{k=1}^N \sum_{i=1}^m \norm{\eikp}^2 + \alpha_2 \sum_{k=1}^N c(k)^2 + \alpha_3, \label{eq:sum_e_bound}
	\end{align}
	where
	\begin{align}
		\alpha_2 &= \frac{2m}{\alpha_1} \left(\frac{m^2 \psi^2}{(1-q)^2} + 4 \right), \nonumber \\
		\alpha_3 &= \frac{\alpha_1}{2} \sum_{i=1}^m \norm{e_i(1)}^2 \nonumber \\
					&+ \frac{2 m^3 \psi^2}{\alpha_1 (1-q)^2}c(0)^2  \nonumber \\
					&+ \frac{2 m \psi q}{1-q} c(1) \sum_{i=1}^m \norm{\li(0)}, \label{eq:constants_lm}
	\end{align}
    and
    \begin{align}
        \psi &= 2 \big ( 1 + \eta^{-(m-1)T} \big ) / \big ( 1 - \eta^{(m-1)T} \big ) \in \bR_+\setminus \{0\} \nonumber \\
        q &= \big ( 1 - \eta^{(m-1)T} \big )^{\frac{1}{(m-1)T}} \in (0,1). \nonumber
    \end{align}
	\begin{proof}
		See the proof of Lemma~3 in \cite{our_paper}.
	\end{proof}
\end{lemma}
It should be noted that for all $i=1,\dots,m$, $\norm{\li(0)}$ is finite as $\li(0)$ is the initialization of the algorithm, $\ai(0)$ is finite since it is the convex combination of finite values, and $x_i(1)$ is finite thanks to Assumption \ref{ass:compact}.
By \eqref{eq:def_error}, and thanks to the fact that $\likp = [\aik + c(k) g_i(\xikp)]_+$ (equation \eqref{eq:subgradient_step}) and to the fact that the projection operator is non-expansive, we have that $\norm{e_i(1)} \leq \norm{c(0)g_i(x_i(1))}$, i.e.  $\norm{e_i(1)}^2$ is finite too.

We then have the following lemma, which is fundamental for the analysis of Section~\ref{sec:proof_of_thm2}.

\begin{lemma} \label{lemma:relation}
	Consider Assumptions~\ref{ass:convex}-\ref{ass:slater} and Assumption \ref{ass:weights}. Fix any $\alpha_1\in\bR_+\setminus\{0\}$. For any $k \in \bN_+$, and for any $x\in X$ and $\lambda\in\bR^p_+$ we have,
		\begin{align}
			\sum_{i=1}^m \norm{&\likp-\lambda}^2 \leq \sum_{i=1}^m \norm{\lik-\lambda}^2 \nonumber \\
				&- (1-\alpha_1)\sum_{i=1}^m \norm{\eikp}^2 \nonumber \\
				&+ \frac{4G^2m}{\alpha_1}c(k)^2 \nonumber \\
				&+ 2 G c(k) \sum_{i=1}^m \norm{\likp - \vkp} \nonumber \\
				&+ 2c(k) \big ( L(x,\vkp) - L(\xkp,\lambda) \big ).\label{eq:general_main_ineq}
		\end{align}
	\begin{proof}
		Consider the quantity $\norm{\aik-\lambda}^2$. Adding and subtracting $\likp$ inside the norm and then expanding the square, we have that\\
		\begin{align}
			\norm{\aik-\lambda}^2 &= \norm{\aik-\likp}^2 + \norm{\likp-\lambda}^2 \nonumber \\
				&+ 2(\aik-\likp)\T (\likp-\lambda) \nonumber \\
				&= \norm{\aik-\likp}^2 + \norm{\likp-\lambda}^2 \nonumber \\
				&+ 2(\aik + c(k)g_i(\xikp) \nonumber \\
					&~~~~~~~~~~~~ - \likp)\T (\likp-\lambda) \nonumber \\
				&- 2c(k)g_i(\xikp)\T (\likp-\lambda), \label{eq:expanding_square}
		\end{align}
		where the second equality is obtained by adding and subtracting $2c(k)g_i(\xikp)\T (\likp-\lambda)$. Consider now step~8 of Algorithm~\ref{alg:Alg1}. By the optimality condition (Proposition 3.1 in \cite[Chapter 3]{bertsekas_tsitsiklis_1997}), we have that, for any $\lambda\in\bR^p_+$,
		\begin{align}
			2(\aik + c(k)g_i(\xikp) &- \likp)\T \nonumber \\
				&\times (\likp-\lambda) \geq 0, \label{eq:optimality_condition_step8}
		\end{align}
		where the first term in the inner product above constitutes the gradient of the objective function that appears at step 8 of Algorithm \ref{alg:Alg1} (it is quadratic, hence differentiable), multiplied by $2c(k)$. Using \eqref{eq:optimality_condition_step8}, we can rewrite \eqref{eq:expanding_square} as an inequality
		\begin{align}
			\norm{\aik-\lambda}^2 &\geq \norm{\aik-\likp}^2 + \norm{\likp-\lambda}^2 \nonumber \\
				&- 2c(k)g_i(\xikp)\T (\likp-\lambda). \label{eq:expanding_square_bis}
		\end{align}
		Now, recalling the definition of $\eikp$, and after rearranging some terms, we have that
		\begin{align}
			\norm{\likp&-\lambda}^2 \leq \norm{\aik-\lambda}^2 - \norm{\eikp}^2 \nonumber \\
				&+ 2c(k)(\likp-\lambda)\T g_i(\xikp),
		\end{align}
		for any $\lambda\in\bR^p_+$. By adding and subtracting $2c(k)(f_i(\xikp)+\aik\T g_i(\xikp))$ in the right-hand side of the inequality above we obtain
		\begin{align}
			\norm{\likp&-\lambda}^2 \leq \norm{\aik-\lambda}^2 - \norm{\eikp}^2 \nonumber \\
				&+ 2c(k)\big( (\likp-\aik)\T g_i(\xikp) \nonumber \\
				&+ f_i(\xikp) + \aik\T g_i(\xikp) \nonumber \\
				&- f_i(\xikp) - \lambda\T g_i(\xikp) \big),\label{eq:lambda_ineq}
		\end{align}
		for any $\lambda\in\bR^p_+$.
		
		Consider now step~7 of Algorithm~\ref{alg:Alg1}. By the optimality of $\xikp$ we have that
		\begin{align}
			f_i(\xikp) + \aik\T g_i&(\xikp) \nonumber \\
&\leq f_i(x_i) + \aik\T g_i(x_i),
		\end{align}
		for any $x_i\in X_i$. Combining the previous statement with \eqref{eq:lambda_ineq}, and by adding and subtracting $2c(k)\vkp\T g_i(x_i)$ and $2c(k)\likp\T g_i(x_i)$, we have that
		\begin{align}
			\norm{\likp&-\lambda}^2 \leq \norm{\aik-\lambda}^2 - \norm{\eikp}^2 \nonumber \\
				&+ 2c(k)\big( (\likp-\aik)\T g_i(\xikp) \nonumber \\
				&+ (\aik - \likp)\T g_i(x_i) \nonumber \\
				&+ (\likp - \vkp)\T g_i(x_i) \nonumber \\
				&+ f_i(x_i) + \vkp\T g_i(x_i)\nonumber  \\
				&- f_i(\xikp) - \lambda\T g_i(\xikp) \big),\label{eq:general_i_ineq}
		\end{align}
		for any $\lambda\in\bR^p_+$ and any $x_i\in X_i$. By summing \eqref{eq:general_i_ineq} across $i$, $i=1,\ldots,m$, rearranging some terms, and recalling the definition of the Lagrangian function and of $\eikp$, we obtain
		\begin{align}
			\sum_{i=1}^m \norm{&\likp-\lambda}^2 \nonumber \\
				&\leq \sum_{i=1}^m \norm{\aik-\lambda}^2 - \sum_{i=1}^m \norm{\eikp}^2 \nonumber \\
				&+ 2c(k) \sum_{i=1}^m \eikp\T (g_i(\xikp)-g_i(x_i)) \nonumber \\
				&+ 2c(k) \sum_{i=1}^m (\likp - \vkp)\T g_i(x_i) \nonumber \\
				&+ 2c(k) \big( L(x,\vkp) - L(\xkp,\lambda) \big), \label{eq:general_ineq}
		\end{align}
		for any $\lambda\in\bR^p_+$ and for any $x\in X$.
		
		By the definition of $\aik$ (step~6 of Algorithm~\ref{alg:Alg1}), by the fact that, under Assumption~\ref{ass:weights}, $\norm{\sum_{j=1}^m a_j^i(k) \ljk -\lambda}^2 = \norm{\sum_{j=1}^m a_j^i(k) (\ljk -\lambda)}^2$, and by convexity of $\norm{\cdot}^2$, we have that
		\begin{equation} \label{eq:general_ineq_term1}
			\sum_{i=1}^m \norm{\aik-\lambda}^2 \leq \sum_{i=1}^m \norm{\lik-\lambda}^2.
		\end{equation}
		Now, from inequality $2a\T b \leq \norm{a}^2 + \norm{b}^2$, where we set $a=\sqrt{\alpha_1}\eikp$ and $b = c(k)(g_i(\xikp)-g_i(x_i))/\sqrt{\alpha_1}$, we obtain
		\begin{align}
			2c(k) &\sum_{i=1}^m \eikp\T (g_i(\xikp)-g_i(x_i)) \nonumber\\
			&\leq \sum_{i=1}^m \alpha_1\norm{\eikp}^2 \nonumber\\
				&~~~~~~~~~~~~~~~~+\sum_{i=1}^m \frac{\norm{g_i(\xikp)-g_i(x_i)}^2}{\alpha_1}c(k)^2  \nonumber\\
			&\leq \sum_{i=1}^m \alpha_1\norm{\eikp}^2 + \sum_{i=1}^m \frac{4G^2}{\alpha_1}c(k)^2 \nonumber\\
			&= \alpha_1 \sum_{i=1}^m \norm{\eikp}^2 + \frac{4G^2m}{\alpha_1}c(k)^2, \label{eq:general_ineq_term2}
		\end{align}
		where the second inequality is given by the fact that $\norm{g_i(\xikp)-g_i(x_i)}\leq 2G$.
		By the Cauchy-Schwarz inequality we have that
		\begin{align}
			2c(k) \sum_{i=1}^m &(\likp - \vkp)\T g_i(x_i) \nonumber  \\
				&\leq 2 G c(k) \sum_{i=1}^m \norm{\likp - \vkp}.\label{eq:general_ineq_term3}
		\end{align}
		Finally, by using \eqref{eq:general_ineq_term1}, \eqref{eq:general_ineq_term2}, \eqref{eq:general_ineq_term3} together with \eqref{eq:general_ineq}, inequality \eqref{eq:general_main_ineq} follows, thus concluding the proof.
	\end{proof}
\end{lemma}

The relations established in Lemmas~\ref{lemma:sum_error} and \ref{lemma:relation} can be exploited to prove the following proposition.

\begin{prop} \label{prop:conv_error}
	Consider Assumptions~\ref{ass:convex}-\ref{ass:network}. We have that
	\begin{enumerate}
		\item $\sum_{k=1}^{\infty} \sum_{i=1}^m \norm{\eik}^2 < \infty$, \\[-0.5em]
		\item $\lim_{k \rightarrow \infty} \norm{\eik} = 0$, for all $i=1,\dots,m$, \\[-0.5em]
		\item $\sum_{k=1}^\infty c(k) \sum_{i=1}^m \norm{\likp - \vkp} < \infty$. 
	\end{enumerate}
	\begin{proof}
		Consider \eqref{eq:general_main_ineq} with $\lambda = \ls$ and $x = \xs$, where $(\xs,\ls)$ is an optimal primal-dual pair. By \eqref{eq:saddle_point} we have that $L(\xs,\vkp) - L(\xkp,\ls) \leq 0$, and hence we can drop this term from the right-hand side of \eqref{eq:general_main_ineq}. Fixing $N\in\bN_+$ and summing across $k$, $k=1,\ldots,N$, we have that
			\begin{align}
				\sum_{k=1}^N \sum_{i=1}^m \norm{&\likp-\ls}^2 \leq \sum_{k=1}^N \sum_{i=1}^m \norm{\lik-\ls}^2 \nonumber \\
					&- (1-\alpha_1)\sum_{k=1}^N \sum_{i=1}^m \norm{\eikp}^2 \nonumber \\
					&+ \frac{4G^2m}{\alpha_1} \sum_{k=1}^N c(k)^2 \nonumber \\
					&+ 2 G \sum_{k=1}^N c(k) \sum_{i=1}^m \norm{\likp - \vkp}.
			\end{align}
		By Lemma~\ref{lemma:sum_error}, after some cancellations, and after neglecting some negative terms on the right-hand side, we obtain
			\begin{align}
				(1-\alpha_1(1&+G))\sum_{k=1}^N \sum_{i=1}^m \norm{\eikp}^2 \nonumber \\
					&\leq \sum_{i=1}^m \norm{\li(1)-\ls}^2 \nonumber \\
					&+ \left( \frac{4G^2m}{\alpha_1} + \alpha_2 G \right) \sum_{k=1}^N c(k)^2 + \alpha_3 G. \label{eq:proof_error}
			\end{align}
		Since \eqref{eq:proof_error} holds for any $\alpha_1 > 0$, 
		one can always choose $\alpha_1$ such that $(1-\alpha_1(1+G)) > 0$.
        Let then $N\rightarrow\infty$. By point (2) of Assumption~\ref{ass:ck_coefficient}, and since $\norm{\li(1)-\ls}^2$ is finite as an effect of $\li(1)$ and $\ls$ being finite (see discussion after the proof of Theorem~\ref{thm:dual_optimality}), the right-hand side of \eqref{eq:proof_error} is finite, leading to point (1) of the proposition. Point (2) then follows directly, while point (3) follows from point (1) together with Lemma~\ref{lemma:sum_error} by letting $N\to\infty$ in \eqref{eq:sum_e_bound}. This concludes the proof.
	\end{proof}
\end{prop}


Interestingly, based on the previous results, we can also prove (see Appendix~\ref{proof_of_summability_dual_objective}) that
\begin{equation} \label{eq:summability_dual_objective}
	\sum_{k=1}^\infty c(k) \sum_{i=1}^m |\varphi_i(\aik)-\varphi_i(\vk)| < \infty.
\end{equation}
Although this result is not necessary to prove Theorem~\ref{thm:primal_optimality}, it is of interest on its own since it allows one to find a lower bound to the rate of convergence to consensus for the dual objective value. As a matter of fact, since
\begin{align}
	&\left(\min_{k=0,\dots,r} \sum_{i=1}^m |\varphi_i(\aik)-\varphi_i(\vk)| \right) \sum_{k=0}^r c(k) \nonumber \\
		&~~~~~~~~~~~~~~~~\leq\sum_{k=1}^r c(k) \sum_{i=1}^m |\varphi_i(\aik)-\varphi_i(\vk)|,
\end{align}
taking the limit for $r\to\infty$ we have that
\begin{equation}
	\lim_{r\to\infty} \left(\min_{k=0,\dots,r} \sum_{i=1}^m |\varphi_i(\aik)-\varphi_i(\vk)| \right) \sum_{k=0}^r c(k) < \infty,
\end{equation}
which implies that the convergence rate of
\begin{equation}
	\min_{k=0,\dots,r} \sum_{i=1}^m |\varphi_i(\aik)-\varphi_i(\vk)|
\end{equation}
to zero cannot be slower than that of $\frac{1}{\sum_{k=0}^r c(k)}$ as $r\!\to\!\infty$.

\subsection{Proof of Theorem~\ref{thm:primal_optimality}} \label{sec:proof_of_thm2}
The proof of Theorem~\ref{thm:primal_optimality} below is inspired by \cite{nedic_ozdaglar_2009_approx}, where the convergence of a running average sequence similar to $\{\hat{x}_i(k)\}_{k \geq 0}$ is studied in a non-distributed setting, though.

Note that $\xhikp$ is a convex combination of past values of $\xikp$, therefore, for all $k\in\bN_+$ we have that $\xhikp\in X_i$.
Consider the quantity $\sum_{i=1}^m g_i(\xhikp)$. By \eqref{eq:xhat_def}, under the convexity requirement of Assumption~\ref{ass:convex}, we have
\begin{align}
	\sum_{i=1}^m g_i(\xhikp) &\leq \sum_{i=1}^m \frac{\sum_{r=0}^k c(r) g_i(\xirp)}{\sum_{r=0}^k c(r)} \nonumber \\
		&= \frac{\sum_{r=0}^k \sum_{i=1}^m c(r) g_i(\xirp)}{\sum_{r=0}^k c(r)}, \label{eq:g_xhat}
\end{align}
where the inequality (as well as the subsequent ones) is to be interpreted component-wise. Step~8 of Algorithm~\ref{alg:Alg1} can be equivalently written as $\likp = [\aik + c(k) g_i(\xikp)]_+$, where $[\,\cdot\,]_+$ denotes the projection of its argument on $\bR^p_+$. (see also discussion at the end of Section \ref{sec:setup}). Therefore,
\begin{equation} \label{eq:lambda_proj_ineq}
	\likp \geq \aik + c(k)g_i(\xikp).
\end{equation}
Summing \eqref{eq:lambda_proj_ineq} with respect to agents and steps, and then substituting in \eqref{eq:g_xhat} gives
\begin{align}
	\sum_{i=1}^m &g_i(\xhikp) \nonumber \\
		&\leq \frac{\sum_{r=0}^k \sum_{i=1}^m (\lirp - \air)}{\sum_{r=0}^k c(r)} \nonumber \\
		&= \frac{\sum_{r=0}^k \sum_{i=1}^m (\lirp - \sum_{j=1}^m a_j^i(k) \ljr)}{\sum_{r=0}^k c(r)} \nonumber \\
		&= \frac{\sum_{r=0}^k \Big ( \sum_{i=1}^m \lirp - \sum_{j=1}^m \sum_{i=1}^m a_j^i(k) \ljr \Big ) }{\sum_{r=0}^k c(r)} \nonumber \\
		&= \frac{\sum_{r=0}^k \Big ( \sum_{i=1}^m \lirp - \sum_{j=1}^m \ljr \Big )}{\sum_{r=0}^k c(r)} \nonumber \\
		&= \frac{\sum_{i=1}^m (\likp - \li(0))}{\sum_{r=0}^k c(r)}, \label{eq:g_xhat_ineq}
\end{align}
where the first equality follows from the definition of $\air$, the second inequality involves an exchange on the summation order, the third equality is due to Assumption~\ref{ass:weights}, and the last one is obtained after some term cancellations. Since $\{\lik\}_{k\geq 0}$ is a bounded sequence (see the discussion after the proof of Theorem~\ref{thm:dual_optimality}), and due to the fact that $\sum_{r=0}^\infty c(r) = \infty$, taking the limit superior in \eqref{eq:g_xhat_ineq} we obtain that
\begin{equation} \label{eq:asymp_feasib}
	\limsup_{k\rightarrow\infty} \sum_{i=1}^m g_i(\xhikp) \leq 0.
\end{equation}
	
Consider now the quantity $2\sum_{i=1}^m L_i(\xhikp,\ls)$ for any $\ls\in\Ls$. By \eqref{eq:xhat_def}, under Assumption~\ref{ass:convex}, we have that
\begin{align}
	2\sum_{i=1}^m L_i(\xhikp,\ls) &\leq 2\sum_{i=1}^m \frac{\sum_{r=0}^k c(r) L_i(\xirp,\ls)}{\sum_{r=0}^k c(r)} \nonumber \\
			&= \frac{\sum_{r=0}^k 2 c(r) L(\xrp,\ls)}{\sum_{r=0}^k c(r)}. \label{eq:L_xhat}
\end{align}
By \eqref{eq:general_main_ineq} in Lemma~\ref{lemma:relation} with $x = \xs$ and $\lambda = \ls$, for any $(\xs,\ls)\in\Xs\x\Ls$, with $r$ in place of $k$, and after neglecting the negative term $-(1-\alpha_1)\sum_{i=1}^m \norm{\eikp}^2$, we have that
\begin{align}
	2 c(r) L(&\xrp,\ls) \nonumber \\
		\leq 2& c(r) L(\xs,\vrp) \nonumber \\
		&+ \sum_{i=1}^m \norm{\lir-\ls}^2 - \sum_{i=1}^m \norm{\lirp-\ls}^2 \nonumber \\
		&+ \frac{4G^2m}{\alpha_1} c(r)^2 \nonumber \\
		&+ 2 G c(r) \sum_{i=1}^m \norm{\lirp - \vrp} \nonumber \\
		\leq 2& c(r) L(\xs,\ls) \nonumber \\
		&+ \sum_{i=1}^m \norm{\lir-\ls}^2 - \sum_{i=1}^m \norm{\lirp-\ls}^2 \nonumber \\
		&+ \frac{4G^2m}{\alpha_1} c(r)^2 \nonumber \\
		&+ 2 G c(r) \sum_{i=1}^m \norm{\lirp - \vrp}, \label{eq:optimal_main_ineq}
\end{align}
where the second inequality follows from the fact that $L(\xs,\vrp) \leq L(\xs,\ls)$ due to \eqref{eq:saddle_point}. Substituting \eqref{eq:optimal_main_ineq} in \eqref{eq:L_xhat} we have that
\begin{align}
	2 L(&\xhkp,\ls) \nonumber \\
		&\leq \frac{\sum_{r=0}^k 2 c(r) L(\xs,\ls)}{\sum_{r=0}^k c(r)} \nonumber \\
		&+ \frac{1}{\sum_{r=0}^k c(r)} \Bigg( \frac{4G^2m}{\alpha_1} \sum_{r=0}^k c(r)^2 \nonumber \\
		&+ \sum_{i=1}^m \norm{\li(0)-\ls}^2-\sum_{i=1}^m \norm{\li(k+1)-\ls}^2 \nonumber \\
		&+ 2 G \sum_{r=0}^k c(r) \sum_{i=1}^m \norm{\lirp - \vrp} \Bigg). \label{eq:L_xhat_ineq}
\end{align}
Using Assumption~\ref{ass:ck_coefficient} point (2), the boundedness of $\{\lik\}_{k\geq 0}$, and Proposition~\ref{prop:conv_error} part (3), we know that all terms inside the parentheses are finite as $k\rightarrow\infty$. Therefore,
\begin{equation}
	\limsup_{k\rightarrow\infty} L(\xhkp,\ls) \leq L(\xs,\ls).
\end{equation}
However, by \eqref{eq:saddle_point} we have that $L(\xhkp,\ls) \geq L(\xs,\ls)$, hence
\begin{equation} \label{eq:asymp_optimality}
	\lim_{k\rightarrow\infty} L(\xhkp,\ls) = L(\xs,\ls).
\end{equation}
	
By \eqref{eq:asymp_feasib} and \eqref{eq:asymp_optimality}, and due to the fact that $L(\cdot,\ls)$ is continuous as an effect of $f_i(\cdot)$ and $g_i(\cdot)$ being convex under Assumption \ref{ass:convex}, we have that all limit points of $\{\xhk\}_{k\geq 0}$ are feasible and achieve the optimal value. This implies that they are optimal for the primal problem, thus concluding the proof.\hfill\ \qed

\section{Numerical Example} \label{sec:example}
In this section we demonstrate the efficacy of the proposed approach on a modified version of the Plug-in Electric Vehicles (PEVs) charging problem described in \cite{vujanic_etal_2016}. This problem consists in finding an optimal overnight charging schedule for a fleet of $m$ vehicles, which has to be compatible with local requirements and limitations (e.g., desired final state of charge and maximum charging power for each vehicle), and must satisfy some network-wide constraints (e.g., maximum power that the network can deliver).

\begin{figure}[b]
	\centering
	\includegraphics[width=0.75\columnwidth]{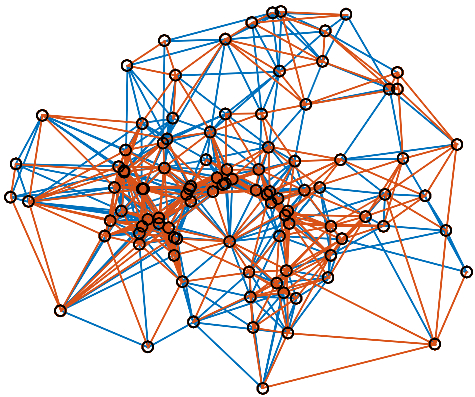}
	\caption{Network of $m=100$ agents.}
	\label{fig:graph}
\end{figure}

Specifically, we hereby consider a slight modification of the ``only charging'' problem in \cite{vujanic_etal_2016}, in that we allow for the optimization of the vehicles charging rate at each time slot, instead of deciding whether to charge or not to charge the vehicle at some fixed charging rate. The overall charging problem can be formalized as the following optimization program
\begin{equation} \label{eq:example_program}
	\begin{aligned}
		\min_{{\{x_i\in X_i\}}_{i=1}^m} \quad &\sum_{i=1}^m c_i\T x_i \\
			\text{subject to:} \quad &\sum_{i=1}^m \left( A_i x_i - \frac{b}{m} \right) \leq 0
	\end{aligned}
\end{equation}
which is a linear program ($X_i$ are indeed bounded convex polytopic sets) having the same structure of \eqref{eq:cP} and satisfying Assumptions~\ref{ass:convex}-\ref{ass:slater}. In \eqref{eq:example_program} the components of the optimization vector $x_i$ represent the charging rate for vehicle $i$ in given time slots, vector $c_i$ gives the costs for charging vehicle $i$ with unitary charging rate, $X_i$ expresses local requirements and limitations for vehicle $i$ such as desired final state of charge and battery rated capacity, while $\sum_{i=1}^m \left( A_i x_i - b/m \right) \leq 0$ encodes network-wide power constraints. We refer the reader to \cite{vujanic_etal_2016} for the precise formulation of all quantities in \eqref{eq:example_program}.

In our simulation we considered a fleet of $m = 100$ vehicles. According to the ``only charging'' set-up in \cite{vujanic_etal_2016}, each vehicles has $n_i = 24$ decision variables and a local constraint set defined by $197$ inequalities. There are $p = 48$ coupling inequalities, and therefore we have $48$ Lagrange multipliers to optimize for the dual problem.

The communication network is depicted in Figure~\ref{fig:graph} and corresponds to a connected graph, whose edges are divided into two groups: the blue and the red ones, which are activated alternatively; this way Assumption~\ref{ass:network} is satisfied with a period of $T = 2$. For each set of edges we created a doubly stochastic matrix so as to satisfy Assumption~\ref{ass:weights}. Finally, we selected $c(k) = \frac{1}{k+1}$.

\begin{figure}[t]
	\centering
	\includegraphics[width=\columnwidth]{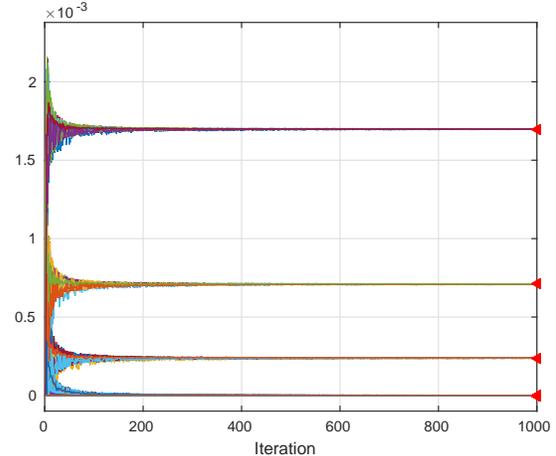}
	\caption{Evolution of the agents' estimates $\lik$, $i=1,\dots,m$. Red triangles represent the optimal dual solution.}
	\label{fig:multipliers}
\end{figure}

\begin{figure}[t]
	\centering
	\includegraphics[width=\columnwidth]{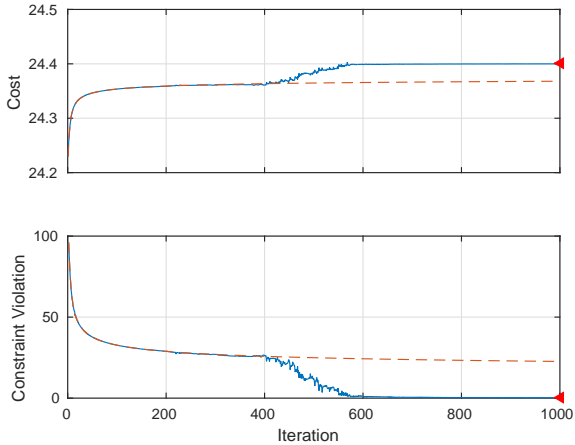}
	\caption{Evolution of primal objective $\sum_{i=1}^m c_i\T x_i$ (upper plot) and constraint violation $\max\{\sum_{i=1}^m (A_i x_i - b/m),0\}$ (lower plot) as a function of $\xh_i(k)$ (dashed lines), and $\tilde{x}_i(k)$ (solid lines).}
	\label{fig:obj_constr}
\end{figure}

We ran Algorithm~\ref{alg:Alg1} for 1000 iterations. Figure~\ref{fig:multipliers} shows the evolution of the agents' estimates $\lik$, $i=1,\dots,m$ across iterations. As expected, all agents gradually reach consensus on the optimal Lagrange multipliers of \eqref{eq:example_program} (red triangles).
Note that, for the problem at hand, only $3$ multipliers are positive, while all the remaining $45$ are equal to zero (in the figure, there are $45$ red triangles in $0$ each one on top of the other). Figure~\ref{fig:obj_constr} instead shows the evolution of the primal objective value $\sum_{i=1}^m c_i\T x_i$ (upper plot), and constraint violation in terms of $\max\{\sum_{i=1}^m (A_i x_i - b/m),0\}$ (lower plot), where $x_i$ is replaced by two different sequences: $\xh_i(k)$ (dashed lines), and $\tilde{x}_i(k)$ (solid lines), $\tilde{x}_i(k)$ being defined as
\begin{equation} 
	\tilde{x}_i(k+1) = \begin{cases}
		\xhikp &k < k_{s,i} \\
		\displaystyle \frac{\sum_{r=k_{s,i}}^k c(r) \xirp}{\sum_{r=k_{s,i}}^k c(r)}  &k \geq k_{s,i}
	\end{cases}
\end{equation}
where $k_{s,i}\in\bN_+$ is the iteration index related to a specific event, namely, the ``practical'' convergence of the Lagrange multipliers, as detected by agent $i$.
Specifically, in the proposed example $k_{s,i}$ is the first iteration step at which the quantity $\norm{\likp-\aik}_2$ has kept below a certain threshold ($10^{-5}$ in our simulation) for $m=100$ consecutive iterations. Being $\tilde{x}_i(k)$ a refresh of $\xh_i(k)$, it is easy to show via the same argument used for $\xh_i(k)$ that Theorem~\ref{thm:primal_optimality} holds also for $\{ \tilde{x}_i(k) \}_{k \geq 0}$, $i=1,\dots,m$.

As can be seen from Figure~\ref{fig:obj_constr}, the rate of convergence of the cost and the constraints violation computed with the $\{\xh_i(k)\}_{k\geq0}$ sequence appears to be slow. In comparison with the rate of $O(1/k)$ in \cite{simonetto2016primal}, a lower bound for our auxiliary sequence to converge is in fact given by $1/\sum_{i=0}^\infty c(k) \sim O(1/\log(k))$ (see the discussion below \eqref{eq:summability_dual_objective}). We therefore believe that the difference in the convergence rate between Algorithm~1 and \cite{simonetto2016primal} might be primarily due to the constant vs. vanishing step-size. Having a vanishing step-size, however, allows us to provide optimality guarantees, while, as discussed in the introduction, in \cite{simonetto2016primal} only convergence to a neighborhood of the optimal solution is guaranteed.
The motivation for introducing the modified auxiliary sequence (which has the same asymptotic convergence rate of the original one) is mainly to counteract the fact that the convergence of $\hat{x}_i(k)$ is also adversely affected by the bad estimates of the Lagrange multipliers obtained at the early stages of the algorithm. By the re-initialization mechanism, $\tilde{x}_i(k)$ for $k\geq k_{s,i}$ depends only on estimates of the Lagrange multipliers that are very close to $\lambda^\star$ and, as such, it presents a much better numerical behavior than $\hat{x}_i(k)$.

\color{black}
\section{Concluding remarks} \label{sec:conclusion}
In this paper we proposed a novel distributed algorithm for a certain class of convex optimization programs, over time-varying multi-agent networks. More precisely, an iterative scheme combining dual decomposition and proximal minimization was conceived, which converges to some optimal dual solution of the centralized problem counterpart, while primal iterates generated by the algorithm converge to the set of primal minimizers. A realistic example on electric vehicles charging over a network was also provided to better illustrate the features of the proposed methodology.

Future work will focus on the analysis of the convergence rate, and on the relaxation of the convexity assumption by extending the results of \cite{vujanic_etal_2016,udell_boyd_2015} to a distributed set-up and quantifying the duality gap incurred in case of mixed-integer programs.

From an application point of view, our goal is to apply the proposed algorithm to the problem of optimal energy management of a building network \cite{Ioli_etal_2015}.

\appendix
\section{Appendix}

\begin{proofof}{\eqref{eq:summability_dual_objective}} \label{proof_of_summability_dual_objective}
	By Theorem~\ref{thm:dual_optimality}, $\{\norm{\lik}\}_{k\geq 0}$ is bounded (see the discussion after the proof of Theorem~\ref{thm:dual_optimality}), whereas by the definition of $\aik$ (step~6 of Algorithm~\ref{alg:Alg1}) and $\vk$ in \eqref{eq:vk}, it follows that $\{\norm{\vk}\}_{k\geq 0}$ and $\{\norm{\aik}\}_{k\geq 0}$ are also bounded, for all $i=1,\dots,m$. Let $\bar{D}\in\bR_+$ denote a uniform upper bound for these sequences. Due to compactness of $X_i$, $i=1,\dots,m$, $\varphi_i(\bar{\lambda})$ is finite for any $\bar{\lambda}\in\bD$, with $\bD = \{\lambda\in\bR^p_+ : \norm{\lambda}\leq \bar{D} \}$. Therefore, $\varphi_i(\cdot)$ is concave (being a dual function) on the compact set $\bD$, hence it will also be Lipschitz continuous on $\bD$ with Lipschitz constant $C_i\in\bR_+$, i.e.,
	\begin{equation}
		|\varphi_i(\lambda_1)-\varphi_i(\lambda_2)| \leq C_i \norm{\lambda_1-\lambda_2}, \; \forall \lambda_1,\lambda_2\in\bD
	\end{equation}
	
	By the definition of $\bar{D}$ we have that $\aik,\vk\in\bD$, for all $k\in\bN_+$, for all $i=1,\dots,m$, hence
	\begin{equation}
		\sum_{i=1}^m |\varphi_i(\aik)-\varphi_i(\vk)| \leq C \sum_{i=1}^m \norm{\aik-\vk},
	\end{equation}
	where $C = \max_{i=1,\dots,m} C_i$. Multiplying both sides by $2c(k)$, fixing $N\in\bN_+$ and summing across $k$, $k=1,\ldots,N$, we have that
	\begin{align}
		2\sum_{k=1}^N c(k) \sum_{i=1}^m &|\varphi_i(\aik)-\varphi_i(\vk)| \nonumber \\
			&\leq 2C\sum_{k=1}^N c(k) \sum_{i=1}^m \norm{\aik-\vk} \nonumber \\
			&\leq 2C\sum_{k=1}^N c(k) \sum_{i=1}^m \norm{\lik-\vk}, \label{eq:dual_sum_relation}
	\end{align}
	where the second inequality follows from
	\begin{align}
		\sum_{i=1}^m \norm{\aik-\vk} &= \sum_{i=1}^m \left\| \sum_{j=1}^m a_j^i(k) (\ljk -\vk) \right\| \nonumber \\
			&\leq \sum_{i=1}^m \sum_{j=1}^m a_j^i(k) \norm{\ljk-\vk} \nonumber \\
			&= \sum_{j=1}^m \norm{\ljk-\vk},
	\end{align}
	where the first equality is obtained by the definition of $\aik$ (step~6 of Algorithm~\ref{alg:Alg1}) and by Assumption~\ref{ass:weights}, the inequality is due to the triangle inequality for $\norm{\cdot}$, and the last equality is obtained exchanging the two summations and using Assumption~\ref{ass:weights}.
	Letting $N\rightarrow\infty$ in \eqref{eq:dual_sum_relation}, and due to Proposition~\ref{prop:conv_error} part~3, \eqref{eq:summability_dual_objective} follows, thus concluding the proof.
\end{proofof}

\bibliographystyle{abbrv}

\end{document}